\documentclass[12pt]{article}
\usepackage{amsmath,amsthm,amsfonts}
\theoremstyle{plain}
\newtheorem{theor}{Theorem}[section]
\theoremstyle{plain}
\newtheorem{lemma}[theor]{Lemma}
\theoremstyle{plain}
\newtheorem{cor}[theor]{Corollary}
\theoremstyle{remark}
\newtheorem{rem}[theor]{Remark}
\newtheorem{example}[theor]{Example}
\theoremstyle{remark}
\newtheorem{defin}[theor]{Definition}
\theoremstyle{definition}
\normalbaselines
\begin{document}
\vbox {\vspace{6mm}}
\begin{center}
{\large \bf Selection of subsystems of random variables\\
equivalent in distribution to the Rademacher system\\[3mm]}
{\bf S.V.Astashkin\\[5mm]}
\end{center}
\begin{abstract}
We present necessary and sufficient conditions on systems of random
variables for them to possess a lacunary subsystem equivalent
in distribution to the Rademacher system on the segment
$[0,1].$
In particular, every uniformly bounded orthonormal
system has this property. Furthermore, an arbitrary finite uniformly
bounded orthonormal set of functions
$\{f_n\}_{n=1}^N$
contains a subset of "logarithmic" density equivalent in distribution to the
corresponding set of Rademacher functions, with a constant independent of
$N.$
A connection between the tail distribution and the
$L_p$
-norms of polynomials with respect to systems of random variables exploited.
We use, also, these results to study the
${\cal K}$
-closed representability of some Banach couples.
\end{abstract}
\def\ab{(\Omega,\Sigma,{\mathbb P})}
\def\bc{\{f_n\}_{n=1}^{\infty}}
\def\cd{\{f_{n_k}\}}
\def\de{\omega\in\Omega}
\def\ef{\sum_{n=1}^{m}a_nf_{n}}
\def\fg{{\mathbb N}}
\def\gh{{\cal K}}
\def\hj{{\mathbb R}}
\def\ji{{\cal J}}
\def\ik{(a_n)_{n=1}^{\infty}}
\def\kl{\{A_j\}_{j=1}^t}
\def\lm{\{r_n\}_{n=1}^\infty}
\def\mn{\{f_n\}}
\def\no{\{f_{n_k}\}}
\def\op{\sum_{n=1}^\infty}
\def\pq{\{r_n\}}
\def\qr{\{\varphi_i\}_{i=1}^{\infty}}
\def\rs{\sum_{n=1}^ma_nr_n}
\def\st{\kappa(t,a)}
\def\tu{|f(\omega)|}
\def\uv{\|a\|_1}
\def\vw{|f_n(\omega)|}
\def\wx{\sum_{i=1}^\infty a_i\varphi_i}
\def\xy{{\cal K}_{1,2}(\sqrt{t},a)}
\def\yz{(a_i)_{i=1}^\infty}
\def\za{\sum_{i=1}^\infty}
\def\yb{\varepsilon_i}
\def\xc{\{\varphi_i\}}
\def\xz{\int_{\Omega}}
\def\ve{\alpha_{i,N}}
\def\uf{\gamma_{i,N}}
\def\tg{\sum_{i\in A_j^N}}
\def\sh{\sum_{j=1}^t}
\def\rj{\sum_{i=N+1}^\infty}
\def\qi{\prod_{i=1}^N}
\def\pk{\prod_{i\in A_j^N}}
\def\oa{\varphi_i(\omega)}
\def\ob{R_N(\omega)}
\def\oc{\|\varphi\|_t}
\def\od{\sum_{i=1}^m}
\def\vr{{\mathbb P}}
\def\xs{{\mathbb E}}

\def\abc{\sum_{i=1}^s}
\def\abd{f_{n_i}}
\def\abe{\{f_{n_i}\}_{i=1}^s}
\def\abf{\{f_n\}_{n=1}^N}
\def\abg{{{\cal A}_s}}
\def\abh{\sum_{\theta\in{\cal A}_s'}}
\def\cdh{\sum_{\theta\in{\cal A}_s}}
\def\abi{{\theta_i}}
\def\abk{\prod_{i=1}^s}
\def\abl{{\Delta_k}}
\def\abm{{\bar g}}
\def\abn{{\alpha_k}}
\def\abo{{\tilde h}}
\def\abp{{[a_{k-1},\alpha_k)}}
\def\abq{{[\alpha_k,a_k]}}
\def\abr{(a_i)_{i=1}^s}
\def\abj{\{h_i\}_{i=1}^s}
\def\abt{\sum_{i\in A_j}}
\def\abu{f_{n_k}}
\def\abv{\prod_{k=1}^s}
\def\abw{{\theta_k}}
\def\abx{\gamma_i}

{\bf 1.Introduction.}
\setcounter{section}{1}
In this paper, we shall consider problems connected to the selection
of lacunary subsystems of random variables.
\begin{defin}
A sequence of random variables (r.v.)
$\bc,$
$f_n\in L_p$
$(p>2)$
defined on a probability space
$\ab$
is said to be a
$\Lambda(p)$
-system
if there exists a constant
$K_p>0$
such that
$$
\Bigl\|\ef\Bigr\|_p\;\le\;K_p\Bigl\|\ef\Bigr\|_2,\eqno{(1.1)}$$
for all
$m\in\fg$
and
$a_n\in\hj$
$(n=1,2,.,m).$
\end{defin}

A sequence
$\bc$
is called a
$\Lambda(\infty)$
-system if it is a
$\Lambda(p)$
-system for any
$p<\infty.$

Another lacunary condition is connected to the absolute convergence of
series of r.v.
\begin{defin}
A system
$\bc$
of r.v. is said to be a Sidon system if
$$
\sum_{n=1}^m|a_n|\;\le\;C\,\Bigl\|\ef\Bigr\|_\infty,\eqno{(1.2)}$$
where a constant
$C>0$
is independent of
$m\in\fg$
and
$a_n\in\hj$
$(n=1,2,.,m).$
\end{defin}

Now we introduce the following concept.
\begin{defin}
We shall say that two systems
$\bc$
and
$\{g_n\}_{n=1}^{\infty}$
of r.v. defined on probability spaces
$\ab$
and
$(\Omega',\Sigma',\vr'),$
respectively, are equivalent in distribution (write:
$\mn\stackrel{\vr}{\sim}\{g_n\}$)
if there exists a constant
$C>0$
such that
$$
C^{-1}\vr\biggl\{\biggl|\sum_{n=1}^{m} a_n f_n(\omega)\biggr|>Cz\biggr\}\;\le\;\vr'\biggl\{\biggl|\sum_{n=1}^{m} a_n g_n(\omega')\biggr|>z\biggr\}\;\le$$
$$
\le\;C \vr\biggl\{\biggl|\sum_{n=1}^{m} a_n f_n(\omega)\biggr|>C^{-1}z\biggr\},\eqno{(1.3)}$$
for all
$m\in\fg,$
$a_n\in\hj$
$(n=1,2,.,m),$
and
$z>0.$
Every possible value of
$C$
will be called a constant of this equivalence.
\end{defin}

The main aim of this paper is to find necessary and sufficient conditions on
systems of r.v. for them to possess a lacunary subsystem equivalent
in distribution to the Rademacher system
$\pq_{n=1}^\infty,$
$$
r_n(x)=\;\rm sign\,\sin(2^{n-1}\pi x),\;\;x\in [0,1].$$
In particulary, we shall show that every uniformly bounded orthonormal
system has this property. This strengthens the following theorem from
the paper [1].
\begin{theor}
Suppose
$\bc$
is an orthonormal system of r.v. on a probability space
$\ab,$
$|f_n(\omega)|\le M,$
$\de,$
$n\in\fg.$
Then there exist a subsystem
$\cd\subset\bc$
and a positive constant
$C=C(M)$
such that for all
$m\in\fg,$
$a_k\in\hj$
$(k=1,2,.,m),$
and
$z>0$
we have
$$
\vr\biggl\{\biggl|\sum_{k=1}^{m} a_k f_{n_k}(\omega)\biggr|>z\biggr\}\;\le\;C\biggl|\biggl\{x\in [0,1]:\,\,\biggl|\sum_{k=1}^{m} a_k r_k(x)\biggr|>C^{-1}z\biggr\}\biggr|$$
(here
$|E|$
denotes Lebesgue measure of a set
$E$).
\end{theor}
Khintchine's inequality implies that the Rademacher system is a
$\Lambda(\infty)$
-system. It easy to check that this system also is a Sidon system.
Therefore a system of r.v. equivalent in distribution to
the Rademacher system is both a
$\Lambda(\infty)$
-system (futhermore, a constant
$K_p$
in (1.1) can be taken the same, up to order, as that for the
Rademacher system, i.e.,
$K_p\,\asymp\,\sqrt{p}$)
and a Sidon system. So, the results of this paper strengthen
well-known Banach's and Gaposhkin's theorems about the selection of
$\Lambda(p)-$
and Sidon subsystems, respectively ([2], [3,Th.1.4.1]).

Similar problems for finite sets of r.v. also are of interest.
Let
$\cal M$
be some class of sets
$\abf$
of r.v., with fixed
$N\in \fg.$
For some lacunary property, we seek as large as possible
$m=m(N)\in\fg $
such that for every system from
$\cal M$
there exists a subset
$S\subset \{1,2,..,N\}$
having the following properties: 1) the number of its elements is not less than
$m;$
2) the subsystem
$\{f_n\}_{n\in S}$
satisfies considered property, with a constant independent of
$N.$

We shall consider this problem for the class of orthonormal sets of functions
defined on the segment
$[0,1]$
and the property of equivalence in distribution to the Rademacher system.
We prove that every such set contains a subset of "logarithmic"
density
$m\asymp\log_2N$
that is equivalent in distribution to the set of first
$m$
Rademacher functions. This result is exact, up to order. It strengthens
Kashin's theorem about the selection of Sidon subsystems (see the
paper [4]).

Let us remark, also, that lacunary subsystems of r.v. equivalent to the
Rademacher system are "best possible". The reason for this is that
the distribution of polynomials corresponding to any subsystem of the
Rademacher system is the same as one for itself.

The important tools in following are Peetre 's
$\gh$
-functional and the real interpolation method of operators.

Let
$(X_0,X_1)$
be a Banach couple,
$x\in X_0+X_1,$
$t>0.$
The Peetre
$\gh$
-functional concerning to this Banach couple is defined by the formula
$$
\gh(t,x;X_0,X_1)=\;\inf\{\|x_0\|_{X_0}\,+\,t\|x_1\|_{X_1}:\:x=x_0+x_1,x_0\in X_0,x_1\in X_1\}.$$
For any
$0<\theta<1,$
$1\le q<\infty$
define the norm
$$
\|x\|_{\theta,q}=\;\left\{\sum_{n=-\infty}^\infty\Bigl[2^{-n\theta}\gh(2^n,x;X_0,X_1)\Bigr]^q\right\}^{1/q}\;\;(x\in X_0+X_1).$$
Then
$(X_0,X_1)_{\theta,q}=\,\{x\in X_0+X_1:\,\,\|x\|_{\theta,q}<\infty\}$
are interpolation Banach spaces concerning to a Banach couple
$(X_0,X_1)$
(i.e., every linear operator defined on
$X_0+X_1$
and bounded from
$X_i$
into
$X_i$
$(i=0,1)$
is also bounded as operator from
$(X_0,X_1)_{\theta,q}$
into itself).
They are called the spaces of the real interpolation method.

This paper is organized as follows. In Section 2
we shall obtain necessary and sufficient conditions for an arbitrary system
of r.v. to be equivalent in distribution to the Rademacher system. Also,
we consider there some important partial cases of such systems. Section 3
contains the main results of the paper. We prove, here, theorems on
the selection of subsystems equivalent in distribution to the
Rademacher system. In Section 4 similar problems are considered in
the case of finite orthonormal sets of functions. Section 5 is devoted
to an application of previous results to studying the
$\gh$
-closed representability of some Banach couples. Indeed, it was this
application which led us to the results in Sections 2 and 3.

Finally, we introduce some notation. Throughout this paper, a random
variable (r.v.) will be a measurable function from a probability space
to the real line. Suppose
$p>0,$
$f(\omega)$
is a r.v. on a probability space
$\ab;$
then
$$
\|f\|_p=\;(\xs|f|^p)^{1/p}=\;\left\{\xz|f(\omega)|^p\,d\vr(\omega)\right\}^{1/p}.$$
If
$a=\ik$
is a sequence of real numbers, then
$$
\|a\|_p=\;\left\{\op|a_n|^p\right\}^{1/p}.$$

As usual,
$L_p$
consists of all r.v.
$f$
such that
$\|f\|_p<\infty$
and
$l_p$
consists of all sequences
$a=\ik$
such that
$\|a\|_p<\infty.$

If
$A$
and
$B$
are two quantities (that may depend upon certain parameters), we shall write
$A\asymp B$
to mean that there exists a positive constant
$C$
such that
$C^{-1}A\,\le\,B\,\le\,CA.$
By
$|E|$
will be denoted Lebesgue measure of a set
$E\subset\hj.$
\vskip 0.6cm

{\bf 2. Systems of r.v. that are equivalent in distribution
to the Rademacher system.}
\setcounter{section}{2}
\setcounter{theor}{0}
The main result of this section establishes the necessary and sufficient
conditions for an arbitrary system
of r.v. to be equivalent in distribution to the Rademacher system.

Denote by
$\gh_{1,2}(t,a)$
the
$\gh$
-functional
$\gh(t,a;l_1,l_2)$
concerning to the Banach couple
$(l_1,l_2).$
By approximate Holmstedt's formula (see [5,Th.4.1] or [6,${\cal x}\,5.7$]),
there exists a constant
$\alpha>0$
such that
for all
$a=\ik\in l_2$
and
$t>0$
$$
\alpha^{-1}\left\{\sum_{i=1}^{[t^2]}a_i^*\,+\,t\biggl[\sum_{i=[t^2]+1}^\infty(a_i^*)^2\biggr]^{1/2}\right\}\;\le\;\gh_{1,2}(t,a)\;\le$$
$$
\le\;\left\{\sum_{i=1}^{[t^2]}a_i^*\,+\,t\biggl[\sum_{i=[t^2]+1}^\infty(a_i^*)^2\biggr]^{1/2}\right\},\eqno{(2.1)}$$
where
$(a_i^*)_{i=1}^\infty$
is the decreasing rearrangement of the sequence
$(|a_n|)_{n=1}^\infty,$
and
$[z]$
is the largest integer that does not exceed
$z.$

Clearly,
$\st=\,\xy$
is a continuous and increasing function for
$t\ge 0.$
Moreover, formula (2.1) shows that
$$
\lim_{t\to 0+}\st\,=\,0 {\mbox \;\;\; and \;\;\; } \alpha^{-1}\|a\|_1\,\le\,\lim_{t\to +\infty}\st\,\le\,\|a\|_1, \eqno{(2.2)}$$
for any sequence
$a=\ik\in l_2.$

\begin{theor}
Let
$\bc$
be a sequence of r.v. defined on a probability space
$\ab.$
The following conditions are equivalent:

1) $\mn\stackrel{\vr}{\sim}\pq$
($\pq_{n=1}^\infty$
is the Rademacher system on
$[0,1]);$

2)
$$
\Bigl\|\ef\Bigr\|_t\;\asymp\;\Bigl\|\rs\Bigr\|_t,$$
with a constant independent of
$t\in [1,\infty),$
$m\in\fg,$
and
$a_n\in\hj$
$(n=1,.,m);$

3)
$$
\Bigl\|\op a_nf_n\Bigr\|_t\;\asymp\;\st,\eqno{(2.3)}$$
with a constant independent of
$t\in [1,\infty)$
and
$a=\ik\in l_2.$
\end{theor}
\begin{proof}
In the paper [7], P.Hitzhenko proved the equivalence similar to (2.3) for Rademacher
functions. Namely,
$$
\Bigl\|\op a_nr_n\Bigr\|_t\;\asymp\;\st,\eqno{(2.4)}$$
with a constant independent of
$t\in [1,\infty)$
and
$a=\ik\in l_2.$
This implies the equivalence
$2)\Leftrightarrow 3).$
Next, by the
definition of the equivalence in distribution, it is obvious that
$1)\Rightarrow 2).$
Therefore we must
prove only that
$3)\Rightarrow 1).$

Fix
$a=(a_n)_{n=1}^m,$
$m\in\fg$
(we can assume that not all
$a_n$
are equal to zero). We shall show that the tail distribution of a polynomial
$$
f(\omega)=\;\ef(\omega)\;\;(\de)$$
depends only on the norms
$\|f\|_t,$
or, by hypothesis, on the function
$\st=\,\xy$
$(t\ge 1).$

Let
$\beta\ge 1$
be a constant of equivalence (2.3). The
$\gh$
-functional
$\st$
is a concave function with respect to
$t$
[6,${\cal x}\,3.1$].
Therefore from the Paley-Zygmund inequality
[8,${\cal x}\,1.6$] we get for
$t\ge 1$
\begin{multline*}
\vr\{|f(\omega)|\ge(2\beta)^{-1}\st\}\;\ge\;\vr\{|f(\omega)|^t\ge 2^{-t}\|f\|_t^t\}\;\ge\\
\ge\;(1-2^{-t})^2\,\frac{\|f\|_t^{2t}}{\|f\|_{2t}^{2t}}\;\ge\;\beta^{-4t}(1-2^{-t})^2\,\left\{\frac{\st}{\kappa(2t,a)}\right\}^{2t}\;\ge\;(2\beta)^{-4t}.
\end{multline*}
Since
$\st$
increases, then
$$
\vr\{|f(\omega)|\ge (2\beta)^{-1}\st\}\;\ge\;(2\beta)^{-4t-4},\eqno{(2.5)}$$
for all
$t>0.$

On the other hand, by Chebychev's inequality and hypothesis, we have
$$
\vr\{|f(\omega)|\ge \lambda\st\}\;\le\;\lambda^{-t}\left[\frac{\|f\|_t}{\st}\right]^t\;\le\;\left(\frac{\beta}{\lambda}\right)^t,$$
for
$t\ge 1$
and
$\lambda>0.$
In particular, if
$\lambda=\,\beta e,$
then
$$
\vr\{|f(\omega)|\ge \beta e\st\}\;\le\;e^{-t}\;\;\;\;\mbox{for}\;\;\;t\ge 1$$
and
$$
\vr\{|f(\omega)|\ge \beta e\st\}\;\le\;e^{1-t}\;\;\;\;\mbox{for}\;\;\;t>0.\eqno{(2.6)}$$

Following [9], we define the functionals
$$
F(s)=\;\sup\{t>0:\,\st\le s\}\;\;\;\;\mbox{and}\;\;\;\;G(s)=\;\inf\{t>0:\,\st\ge s\}.$$

Suppose that
$0<z<(4\alpha\beta)^{-1}\uv.$
By (2.2),
$\st\ge 4z\beta,$
for some
$t>0.$
Therefore from (2.5) it follows
$$
\vr\{|f(\omega)|>z\}\;\ge\;\vr\{|f(\omega)|\ge (2\beta)^{-1}\st\}\;\ge\;(2\beta)^{-4t-4}.$$
Hence, by the definition of
$G,$
we have
$$
\vr\{|f(\omega)|>z\}\;\ge\;(2\beta)^{-4G(4\beta z)-4}$$
and
$$
\vr\{|f(\omega)|>z\}\;\ge\;C_1^{-1}e^{-C_2G(4\beta z)}\;\;\;\mbox{for}\;\;\;z<(4\alpha\beta)^{-1}\uv,\eqno{(2.7)}$$
where
$C_1=\,(2\beta)^4,$
$C_2=\,4\ln(2\beta).$

Next, as
$\kappa(t,a)$
is a concave function, the definition of
$G$
and
$F$
implies
$$
\kappa(C_2G(4\beta z),a)\;\le\sqrt{2}\,C_2\kappa(2^{-1}G(4\beta z),a)\;\le\;4\sqrt{2}\,C_2\beta z\,=\,C_3z$$
and
$$
C_2G(4\beta z)\;\le\;F(C_3z)\;\;\;\;\mbox{for}\;\;\;z>0,$$
where
$C_3=\,4\sqrt{2}\,C_2\beta.$
If we combine this with (2.7), then we get
$$
\vr\{|f(\omega)|>z\}\;\ge\;C_1^{-1}e^{-F(C_3z)}\;\;\;\;\mbox{for}\;\;\;z<(4\alpha\beta)^{-1}\uv.\eqno{(2.8)}$$

In addition, from (2.3) it follows
$$
\|f_n\|_\infty=\lim_{t\to +\infty}\|f_n\|_t\le\beta,$$
and therefore
$\|f\|_\infty\le\beta\uv.$
Thus we have
$$
\vr\{|f(\omega)|>z\}\;=\;0\;\;\;\;\mbox{for}\;\;\;z\ge\beta\uv.\eqno{(2.9)}$$

Next, by (2.2) we can find
$t'>0$
such that
$0<\kappa(t',a)\le (2\beta e)^{-1}z.$
Then from (2.6) we obtain
$$
\vr\{|f(\omega)|>z\}\;\le\;\vr\{|f(\omega)|\ge\beta e\kappa(t',a)\}\;\le\;e^{1-t'}.$$
The last yields
$$
\vr\{|f(\omega)|>z\}\;\le\;e^{1-F(C_4^{-1}z)}\;\;(z>0),\eqno{(2.10)}$$
where
$C_4=\,2\beta e.$

In the same way, using (2.4), we can prove also the formulas similar
to (2.8) --- (2.10) for the Rademacher functions.
Namely, if
$\beta'\ge 1$
is a constant of equivalence (2.4) and
$r(x)=\,\rs(x),$
then
$$
|\{|r(x)|>z\}|\;\ge\;(C_1')^{-1}e^{-F(C_3'z)}\;\;\;\;\mbox{for}\;\;\;z<(4\alpha\beta')^{-1}\uv,\eqno{(2.8')}$$
$$
|\{|r(x)|>z\}|\;=\;0\;\;\;\;\mbox{for}\;\;\;z\ge\uv,\eqno{(2.9')}$$
and
$$
|\{|r(x)|>z\}|\;\le\;e^{1-F((C_4')^{-1}z)}\;\;\;\;\mbox{for}\;\;\;z>0.\eqno{(2.10')}$$

Let us denote
$$
A=\;\max(\,C_1e\,,\,C_3C_4'\,,\,C_1'e\,,\,C_3'C_4\,,\,4\alpha\beta\beta'\,).\eqno{(2.11)}$$
Note that
$\alpha$
and
$\beta'$
are universal constants. Hence
$A$
depends only on the constant
$\beta$
of equivalence (2.3).

In the case
$z<(4\alpha\beta)^{-1}\uv$
formulas $(2.10'),$ $(2.8),$ and (2.11) give
$$
\{|r(x)|>Az\}|\;\le\;e^{1-F((C_4')^{-1}Az)}\;\le\;e^{1-F(C_3z)}\;\le$$
$$
\le\;C_1e\,\vr\{|f(\omega)|>z\}\;\le\;A\,\vr\{|f(\omega)|>z\}.\eqno{(2.12)}$$
If
$z\ge(4\alpha\beta)^{-1}\uv,$
then
$Az\ge\uv,$
by (2.11). Hence $(2.9')$ yields:
$|\{|r(x)|>Az\}|\,=\,0,$
and inequality (2.12) holds for all
$z>0.$

Conversely, if
$z/A<(4\alpha\beta')^{-1}\uv,$
then from (2.10), $(2.8'),$ and (2.11) we have
$$
\vr\{|f(\omega)|>z\}\;\le\;e^{1-F(C_4^{-1}z)}\;\le\;e^{1-F(C_3'(z/A))}\;\le$$
$$
\le\;C_1'e\,|\{|r(x)|>z/A\}|\;\le\;A\,|\{|r(x)|>z/A\}|.\eqno{(2.13)}$$
If
$z/A\ge(4\alpha\beta')^{-1}\uv,$
then (2.11) implies
$z\ge\beta\uv.$
Taking into account (2.9), we get
$\vr\{|f(\omega)|>z\}\,=\,0,$
and inequality (2.13) is valid again.

From (2.12) and (2.13) it follows that
$\mn\stackrel{\vr}{\sim}\pq.$
This completes the proof of Theorem 2.1.
\end{proof}
\begin{rem}
The proof of Theorem 2.1 shows that constants of the equivalences of
relations 1) --- 3) depend only from each other. Therefore,
for example, if a finite set
$\abf$
satisfies (2.3) with a constant independent of
$N=1,2,..,$
then
$\abf\stackrel{\vr}{\sim}\{r_n\}_{n=1}^N$
also with a constant independent of
$N=1,2,..$
\end{rem}
\vskip 0.3cm

Now, we consider two following partial cases.
\vskip 0.15cm

\begin{defin}
A system of r.v.
$\bc$
 is called multiplicative if for any pairwise distinct
$n_1,n_2,.,n_k$
$(k\in\fg)$
$$
\xs(f_{n_1}f_{n_2}...f_{n_k})\;=\;0.$$
If, in addition,
for any pairwise distinct
$n_1,n_2,.,n_k$
and
$n\ne n_s$
$(s=1,2,.,k)$
$$
\xs(f_{n_1}f_{n_2}...f_{n_k}\,f_n^2)\;=\;0,$$
then
$\mn$
is called a strongly multiplicative system.
\end{defin}

Let us present examples of multiplicative and strongly multiplicative
systems from the paper [10].

\begin{example}
Let
$f_n(x)=\,\sin(2\pi k_nx)$
$(x\in [0,1]).$
If
$k_{n+1}/k_n\ge 2,$
then this sequence is a multiplicative system; moreover, it is a strongly
multiplicative system provided
$k_{n+1}/k_n\ge 3.$
\end{example}
\begin{example}
A sequence
$\bc$
of mean zero and square integrable independent r.v.
is a strongly multiplicative system.
\end{example}
\vskip 0.2cm

In the paper [10], J.Jakubovski and S.Kwapien proved the following theorem.
\begin{theor}[J.Jakubovski and S.Kwapien]
Let
$\{\varphi_n\}_{n=1}^N$
be a multiplicative system and
$\{\psi_n\}_{n=1}^N$
be a strongly multiplicative system of r.v. If
$$
\|\psi_n\|_\infty\,\|\varphi_n\|_\infty\;\le\;\xs(\psi_n^2),\eqno{(2.14)}$$
for each
$n=1,.,N,$
then there exist a probability space
$(\Omega',\Sigma',\vr'),$
$\sigma$
-field
$\Sigma_0\subset\Sigma',$
and a random vector
$(\psi_1',\psi_2',.,\psi_N')$
defined on
$(\Omega',\Sigma',\vr')$
and equidistributed with the vector
$(\psi_1,\psi_2,.,\psi_N)$
such that the random vector
$(\varphi_1,\varphi_2,.,\varphi_N)$
is equidistributed with the conditional expectation
$\xs\Bigl((\psi_1',\psi_2',.,\psi_N')|\Sigma_0\Bigr).$

In particular, if
$H:\;\hj^n\to\hj$
is a convex function, then
$$
\xs\Bigl(H(\varphi_1,.,\varphi_N)\Bigr)\;\le\;\xs\Bigl(H(\psi_1,.,\psi_N)\Bigr).$$
\end{theor}
Theorems 2.6 and 2.1 imply the following statement.
\begin{theor}
Suppose
$\bc$
is a strongly multiplicative system of r.v. Then
$\mn\stackrel{\vr}{\sim}\pq$
(with a constant depending only on
$D$
and
$d$)
iff the following conditions hold:
$$
\vw\le D\;\;\;\;\mbox{for}\;\;\;n=1,2,..\;\;\;{and}\;\;\;{and}\;\;\;\de\eqno{(2.15)}$$
and
$$
d=\,\inf_{n=1,2,..}\xs(f_n^2)>0.\eqno{(2.16)}$$
\end{theor}
\begin{proof}
Suppose that a system
$\{f_n\}_{n=1}^\infty$
satiesfies (2.15) and (2.16). Then for
$\varphi_n=\,f_n/D$
and
$\psi_n=\,r_n$
condition (2.14) is fulfilled. Hence,
$$
\Bigl\|\ef\Bigr\|_t\;\le\;D\,\Bigl\|\rs\Bigr\|_t,$$
for any
$t\ge 1,$
$m\in\fg,$
and
$a_n\in\hj$
$(n=1,2,.,m).$

Analogously, the systems
$\varphi_n=\,(d/D)r_n$
and
$\psi_n=\,f_n$
satisfy (2.14), also. Therefore,
$$
\Bigl\|\rs\Bigr\|_t\;\le\;\frac{D}{d}\,\Bigl\|\ef\Bigr\|_t.$$
Using Theorem 2.1, we obtain:
$\mn\stackrel{\vr}{\sim}\pq,$
and a constant of this equivalnce depends only on
$D$
and
$d$
(see Remark 2.2).

The opposite assertion is an easy consequence of Definition 1.3.
This completes the proof.
\end{proof}
\begin{cor}
A sequence
$\bc$
of mean zero and square integrable independent r.v.
is equivalent in distribution to the Rademacher system if and only if
conditions (2.15) and (2.16) are fulfilled.
\end{cor}

Now consider another situation. Suppose that
$G$
1s a compact abelian group with dual group
$\Gamma.$
Let
$\mu$
be the normalized Haar measure on
$G.$
If
$F\subset\Gamma.$
then by
$C_F(G)$
denote
a set of all continuous functions
$f$
on
$G$
such that
$$
\stackrel{\wedge}{f}(\gamma)=\;\int_G f\stackrel{-}{\gamma}\,d\mu\;=\;0\;\;\;\;\mbox{for}\;\;\;\gamma\not\in F.$$
According to Definition 1.2, a subset
$F$
is called a Sidon set if there exists a constant
$C=C(F)$
(depending only on
$F$)
such that
$$
\sum_{\gamma\in\Gamma}|\stackrel{\wedge}{f}(\gamma)|\;\le\;C\,\|f\|_\infty,$$
for every
$f\in C_F(G).$

In the paper[11], G.Pisier proved the following result.
\vskip 0.2cm
\begin{theor}[G.Pisier]
Suppose that
$F=\{\gamma_n\}\subset\Gamma$
is a Sidon set. Then
$$
\Bigl\|\sum_{n=1}^ma_n\gamma_n\Bigr\|_t\;\asymp\;\Bigl\|\rs\Bigr\|_t\;\;\;\;\mbox{for all}\;\;\;t\ge 1,$$
with a constant depending only on the Sidon constant
$C(F).$
\end{theor}
\vskip 0.2cm

From Theorems 2.9 and 2.1 we obtain the following statement,
which was be proved in a different way in the paper [12] (however,
using the Pisier theorem, also).
\begin{theor}
Every infinite Sidon system
$F=\{\gamma_n\}_{n=1}^\infty$
of characters defined on a compact abelian group is equivalent in distribution
to the Rademacher system on the segment
$[0,1].$
A constant of this equivalence depends only on the Sidon constant
$C(F).$
\end{theor}
\begin{cor}
Sequences
$f_n(x)=\,\sin(2\pi k_nx)$
and
$g_n(x)=\,\cos(2\pi k_nx)$
$(x\in [0,1])$
are equivalent in distribution to the Rademacher system provided
$k_{n+1}/k_n\ge\lambda>1.$
\end{cor}
\begin{rem}
The last result explains why many properties of the Rademacher system
are shared by Hadamard lacunary trigonometric systems.
In particular, all rearrangement invariant norms [see [13] or [14])
of Rademacher
polynomials and lacunary trigonometric polynomials are equivalent
with a certain universal constant. Note that for this reason, for example,
some conditions of theorems of the paper [15] are superflous.
\end{rem}
\vskip 0.6cm

{\bf 3. Selection of subsystems equivalent in distribution to
the Rademacher system.}
\setcounter{section}{3}
\setcounter{theor}{0}
The following approximate formula for the
$\gh$
-functional
$\gh_{1,2}(t,a)$
of S.Montgomery-Smith [16] will serve as a important tool in our proofs.

For arbitrary
$t\in\fg$
define the norm on
$l_2$
$$
\|a\|_{Q(t)}\;=\;\sup\left\{\sh\left(\sum_{n\in A_j}a_n^2\right)^{1/2}\right\},\eqno{(2.1)}$$
where the supremum is taken over all disjoint subsets
$A_1,A_2,..,A_t$
of
$\fg.$
\begin{lemma}[S.Montgomery-Smith]
If
$a=\ik\in l_2$
and
$t^2\in\fg,$
then
$$
\|a\|_{Q(t^2)}\;\le\;\gh_{1,2}(t,a)\;\le\;\sqrt{2}\,\|a\|_{Q(t^2)}.\eqno{(2.2)}$$
\end{lemma}
\begin{theor}
Suppose a sequence of r.v.
$\bc$
defined on a probability space
$\ab$
contains a subsequence
$\cd_{k=1}^\infty$
that satisfies the following conditions:

1) $|f_{n_k}(\omega)|\le D$
$(k=1,2,..;\de);$

2) $f_{n_k}\to 0$
weakly in
$L_2;$

3) $d=\,\inf_{k=1,2,..}\|f_{n_k}\|_2>0.$

Then there exists a subsystem
$\qr\subset\bc$
such that
$$
\Bigl\|\wx\Bigr\|_t\;\asymp\;\xy\;\;\;\;\mbox{for}\;\;\;a=\yz\in l_2\;\;\;\mbox{and}\;\;\;t\ge 1,\eqno{(3.1)},$$
with a constant depending only on
$D$
and
$d.$
\end{theor}
\begin{proof}
Under assumptions 1) and 2) there exists a subsequence
$\{g_i\}\subset\no$
such that
$$
\Bigl\|\za a_ig_i\Bigr\|_t\;\le\;C_1\sqrt{t}\,\|a\|_2,\eqno{(3.2)}$$
where a constant
$C_1=C_1(D)$
does not depend of
$a=\yz\in l_2$
and
$t\ge 1$
([3,Th.1.3.2], detailed proof see in [1]).

Then, if
$(a_i)=\,(b_i)\,+\,(c_i),$
$(b_i)\in l_1,$
and
$(c_i)\in l_2,$
we have
\begin{multline*}
\Bigl\|\za a_ig_i\Bigr\|_t\;\le\;\Bigl\|\za b_ig_i\Bigr\|_\infty\;+\;\Bigl\|\za c_ig_i\Bigr\|_t\;\le\;D\,\za|b_i|\;+\\
+\;C_1\sqrt{t}\,\biggl(\za c_i^2\biggr)^{1/2}\;\le\;\max(C_1,D)\,(\|b\|_1\,+\,\sqrt{t}\|c\|_2).
\end{multline*}
Hence, by the definition of the
$\gh$
-functional,
$$
\Bigl\|\za a_ig_i\Bigr\|_t\;\le\;\max(C_1,D)\,\xy\;\;\;\;\mbox{for}\;\;\;a=\yz\in l_2\;\;\;\mbox{and}\;\;\;t\ge 1.\eqno{(3.3)}$$

In the proof of the opposite inequality we shall use Lemma 3.1 and
an estimate of
$L_q$
-norms
of modified Riesz products from above for
$q>1$
(applications of usual Riesz products to similar problems see, for example,
in [3,${\cal x}\,1.4$] and [17,${\cal x}\,8.4$]).

First note that from conditions 1) and 3) of Theorem 2.1 concerning to
$\cd_{k=1}^\infty$
it follows
$0<d\le \|g_i\|_2\le D$
$(i=1,2,..).$
Therefore, without loss of generality, we can assume that
$\|g_i\|_2\,=\,1$
for all
$i=1,2,..$

Let
$(\yb)_{i=1}^\infty$
be a sequence of real numbers such that
$$
\yb\to 0\;,\;0<\yb<\frac{1}{16}\min(1,D),\;\;\;\mbox{and}\;\;\;\sum_{k=i+1}^\infty\varepsilon_k<\yb\;\;(i=1,2,..).\eqno{(3.4)}$$

It is readily seen that
$\{g_i^2\}$
is a weakly compact sequence in the space
$L_2.$
Hence there are a r.v.
$h=h(\omega)$
and a subsequence
$\{h_k\}\subset\{g_i\}$
such that
$0\le h(\omega)\le D^2,$
$\xs(h)=1,$
and
$h_k^2\to h$
weakly in
$L_2.$
Futhermore,
$h_k\to 0$
weakly in
$L_2.$
Therefore there exists a positive integer
$k_1$
such that
$$
|\xs(h_{k_1})|\;+\;|\xs(h\,h_{k_1})|\;+\;|\xs(h_{k_1}^2-h)|\;\le\;\frac{\varepsilon_1}{2D}.$$
Denote
$\varphi_1=\,h_{k_1}.$
Suppose that positive integers
$k_1<k_2<...<k_{i-1}$
and functions
$\varphi_1=\,h_{k_1},\varphi_2=\,h_{k_2},.,\varphi_{i-1}=\,h_{k_{i-1}}$
$(i\ge 2)$
are chosen. Let
$k_i$
be a positive integer such that
$k_i>k_{i-1}$
and for the function
$\varphi_i=\,h_{k_i}$
$$
\sum\,\biggl\{|\xs(\varphi_{j_1}...\varphi_{j_s}\varphi_i)|\,+|\xs(h\varphi_{j_1}...\varphi_{j_s}\varphi_i)|\,+|\xs[\varphi_{j_1}...\varphi_{j_s}(\varphi_i^2-h)]|\,+$$
$$
+\,\sum_{l=1}^s|\xs(\varphi_{j_1}...\varphi_{j_{l-1}}\varphi_{j_l}^2\varphi_{j_{l+1}}...\varphi_{j_s}\varphi_i)|\biggr\}\;\le\;2^{-i}D^{-1}\yb,\eqno{(3.5)}$$
where
$\varphi_{j_0}=\varphi_{j_{s+1}}=1$
and the summation is taken over all sets of indices
$1\le j_1<j_2<...<j_l<...<j_s\le i-1$
$(s=1,2,.,i-1).$

Let us show that (3.1) is valid for a such sequence
$\qr.$

First let
$t$
be a positive integer. Suppose that
$\kl$
is an arbitrary partition of
$\fg.$
For
$N\in\fg$
and
$j=1,2,.,t$
define sets
$A_j^N=\,\{i=1,.,N:\,i\in A_j\}$
(it is possible some of them are empty). Introduce
block Riesz products corresponding to this partition:
$$
R_N(\omega)=\;\qi(1\,+\,b_i\varphi_i(\omega))\;=\;\prod_{j=1}^t\pk(1\,+\,b_i\varphi_i(\omega)),$$
where
$b_i\in\hj$
such that
$$
\sum_{i\in A_j}b_i^2\;\le\;D^{-2}\eqno{(3.6)}$$

Denote
$$
\varphi(\omega)=\;\wx(\omega)\;\;\;\;\mbox{for}\;\;\;a=\yz\in l_2.$$
Let us estimate the integral
$$
I_N=\;\xz R_N(\omega)\varphi(\omega)\,d\vr(\omega)\;=\;\za a_i\xs(R_N\varphi_i)\;=\;\za a_i\beta_{i,N}\eqno{(3.7)}$$
from below. Here
$$
\beta_{i,N}=\;\xs(R_N\varphi_i)\;=\;\xz\biggl[1\,+\,\sum_{k=1}^N b_k\varphi_k(\omega)\,+\,\sum_{1\le k_1<k_2\le N} b_{k_1}b_{k_2}\varphi_{k_1}(\omega)\varphi_{k_2}(\omega)\,d\vr(\omega)\,+...$$
$$
...+\,b_1b_2...b_N\varphi_1(\omega)...\varphi_N(\omega)\biggr]\varphi_i(\omega)\,d\vr(\omega)\;=\;b_i\ve\;+\;\uf,\eqno{(3.8)}$$
where
$\ve=\,1$
$(i\le N),$
$\ve=\,0$
$(i>N),$
and
\begin{multline*}
\uf=\;\xs(\varphi_i)\;+\;\sum_{k=1,k\ne i}^Nb_k\xs(\varphi_k\varphi_i)\;+\;\sum_{1\le k_1<k_2\le N}b_{k_1}b_{k_2}\xs(\varphi_{k_1}\varphi_{k_2}\varphi_i)\;+\;...\\
...\;+\;\sum_{1\le k_1<...<k_s\le N}b_{k_1}...b_{k_s}\xs(\varphi_{k_1}...\varphi_{k_s}\varphi_i)\;+\;...\;+\;b_1...b_N\xs(\varphi_1...\varphi_N\varphi_i).
\end{multline*}

Then
$$
\uf=\;S_1^i\;+\;S_2^i\;+\;S_3^i,\eqno{(3.9)}$$
where the sums
$S_k^i$
$(k=1,2,3)$
are defined as follows. Note that
$$
\xs(\varphi_{k_1}...\varphi_{k_s}\varphi_i^2)=\;\xs(\varphi_{k_1}...\varphi_{k_s}h)\;+\;\xs[\varphi_{k_1}...\varphi_{k_s}(\varphi_i^2-h)].$$
The sum
$S_1^i$
contains all terms with integrals of the form
$\xs(\varphi_{k_1}...\varphi_{k_s}\varphi_i)$
or
$\xs[\varphi_{k_1}...\varphi_{k_s}(\varphi_i^2-h)],$
where
$k_1<k_2<...<k_s<i.$
In the sums
$S_2^i$
and
$S_3^i$
we include all terms with integrals of the form
$\xs(\varphi_{k_1}...\varphi_{k_{l-1}}\varphi_i\varphi_{k_l}...\varphi_{k_s}),$
($k_1<k_2<...<k_{l-1}\le i<k_l<...<k_s,$
$1\le l\le s$)
and
$\xs(\varphi_{k_1}...\varphi_{k_s}h),$
respectively.

Combining (3.5), (3.4), and (3.6), we get
$$
|S_1^i|\;\le\;\frac{\yb}{D}\;\;,\;\;|S_2^i|\;\le\;\frac{1}{D}\sum_{k=i+1}^\infty\varepsilon_k\;<\;\frac{\yb}{D}.\eqno{(3.10)}$$
Every term of the sum
$S_3^i$
contains the factor
$b_i.$
Therefore from (3.4) we have
$$
|S_3^i|\;\le\;\frac{|b_i|}{D}\za\yb\;\le\;\frac{|b_i|}{8}\eqno{(3.11)}$$
Futhermore, if
$i>N,$
then the sums
$S_2^i$
and
$S_3^i$
in (3.9) are absent. Hence,
$$
\sum_{i=N+1}^\infty |\uf|\;\le\;\frac{1}{D}\sum_{i=N+1}^\infty \yb\;\le\;\frac{1}{16D}.\eqno{(3.12)}$$

Now, let
$1\le i\le N.$
From (3.9) --- (3.11), (3.4), and (3.6)
$$
\biggl(\tg\uf^2\biggr)^{1/2}\,\le\,\sum_{k=1}^3\biggl(\tg(S_k^i)^2\biggr)^{1/2}\,\le\,\frac{2}{D}\tg\yb\,+\,\frac{1}{8}\biggl(\tg b_i^2\biggr)^{1/2}\,\le\,\frac{3}{8D},\eqno{(3.13)}$$
for every
$j=1,2,.,t.$

Formulas (3.7) and (3.8) imply
\begin{multline*}
I_N\;=\;\sum_{i=1}^N a_ib_i\;+\;\sum_{i=1}^Na_i\uf\;+\;\rj a_i\uf\;=\\
=\;\sh\biggl(\tg a_ib_i\biggr)\;+\;\sh\biggl(\tg a_i\uf\biggr)\;+\;\rj a_i\uf.
\end{multline*}
Let
$N\in\fg$
such that
$$
\|a\|_2\;\le\;2\,\biggl(\sum_{i=1}^N a_i^2\biggr)^{1/2}\;\le\;2\,\sum_{j=1}^t\biggl(\tg a_i^2\biggr)^{1/2}.$$
For every
$j=1,2,.,t,$
choose
$b_i$
$(i\in A_j^N)$
such that (3.6) holds and also
$$
\tg a_ib_i\;=\;\frac{1}{D}\biggl(\tg a_i^2\biggr)^{1/2}.$$

By the last equality and inequalities (3.12), (3.13), and (3.4), we obtain
\begin{multline*}
I_N\;\ge\;\frac{1}{D}\sh\biggl(\tg a_i^2\biggr)^{1/2}\;-\;\sh\biggl(\tg a_i^2\biggr)^{1/2}\,\biggl(\tg\uf^2\biggr)^{1/2}\;-\\
-\;\rj|a_i|\,|\uf|\;\ge\;\frac{5}{8D}\,\sh\biggl(\tg a_i^2\biggr)^{1/2}\;-\\
-\;\biggl(\rj|a_i|^2\biggr)^{1/2}\,\biggl(\rj\uf^2\biggr)^{1/2}\;\ge\;\frac{1}{2D}\,\sh\biggl(\tg a_i^2\biggr)^{1/2}.
\end{multline*}
So, for every
$t\in\fg$
and a partition
$\kl$
of positive integers there exist enough large positive integer
$N$
and real numbers
$b_i$
$(i=1,2,.,N)$
satisfying (3.6) such that
$$
I_N=\;\xz R_N(\omega)\varphi(\omega)\,d\vr(\omega)\;\ge\;\frac{1}{3D}\,\sh\left(\sum_{i\in A_j} a_i^2\right)^{1/2}.\eqno{(3.14)}$$

Now, our aim is to prove an estimate of the integral
$I_N$
from above with an expression of the form
$C\|\varphi\|_t,$
with a positive constant
$C.$
By H\"{o}lder's inequality,
$$
|I_N|\;\le\;\|R_N\|_{t'}\,\oc,\eqno{(3.15)}$$
where
$t'=\,t/(t-1).$

Since
$R_N(\omega)\ge 0,$
then it is enough to estimate for this from above the quantity
$$
L_N=\;\|R_N\|_{t'}^{t'}\;=\;\xz[\ob]^{t'}\,d\vr(\omega).\eqno{(3.16)}$$
Let
$t\in\fg,$
$t\ge 3.$
As before,
$\kl$
is a partition of
$\fg,$
$A_j^N=\,\{i=1,.,N:\,i\in A_j\},$
and real numbers
$b_i$
satisfy inequalities (3.6).

From the obvious inequality
$(1+x)^y\,\le\,1+yx\;(x\ge 1,0<y\le 1)$
it follows that
$$
\Bigl[\ob\Bigr]^{t'}\;\le\;\qi(1+b_i\oa)\,[1+(t-1)^{-1}b_i\oa]\;\le$$
$$
\le\;\qi\biggr[1\,+\,\frac{D^2}{t-1}b_i^2\,+\,\frac{t}{t-1}b_i\oa\biggr]\;=\;\prod_{j=1}^t\pk\biggl[1\,+\,\frac{D^2}{t-1}b_i^2\,+\,\frac{t}{t-1}b_i\oa\biggr]\eqno{(3.17)}$$
Denote by
$m(B)$
the number of elements of a set
$B\subset\fg.$
Then the inner product in the last expression is equal to the sum:
\begin{multline*}
\left\{1\,+\,\frac{D^2}{t-1}\tg b_i^2\,+\,\left(\frac{D^2}{t-1}\right)^2\sum_{\stackrel{i_1<i_2}{i_1,i_2\in A_j^N}}b_{i_1}^2b_{i_2}^2\,+...+\,\left(\frac{D^2}{t-1}\right)^{m(A_j^N)}\pk b_i^2\right\}\,+\\
+\,\sum_{C_j\subset A_j^N}\,\left[1\,+\,\left(\frac{D^2}{t-1}\right)^{m(A_j^N)-m(C_j)}\prod_{i\in A_j^N\setminus C_j}b_i^2\right]\,\left(\frac{t}{t-1}\right)^{m(C_j)}\prod_{i\in C_j}b_i\oa,
\end{multline*}
where the summation in the last term is taken over all non-empty subsets
$C_j$
of the set
$A_j^N.$
In view of (3.6) the expression in braces is not greater than
$(t-1)/(t-2),$
if
$t\ge 3.$

So, combining (3.16),(3.17), and (3.6) again, we obtain:
$$
L_N\;\le\;\left(\frac{t-1}{t-2}\right)^t\;+\;\left(\frac{t-1}{t-2}\right)^{t-1}\,\sh\sum_{C_j\subset A_j^N}\,\left(\frac{t}{t-1}\right)^{m(C_j)+1}\left|\xs\left[\prod_{i\in C_j}\oa\right]\right|\;+$$
$$
+\;\left(\frac{t-1}{t-2}\right)^{t-2}\,\sum_{1\le j_1<j_2\le t}\sum_{C_{j_1}\subset A_{j_1}^N}\sum_{C_{j_2}\subset A_{j_2}^N}\left(\frac{t}{t-1}\right)^{m(C_{j_1})+m(C_{j_2})+2}\,\left|\xs\left[\prod_{k=1}^2\prod_{i\in C_{j_k}}\oa\right]\right|\;+...$$
$$
..+\;\left(\frac{t-1}{t-2}\right)^{t-s}\,\sum_{1\le j_1<...<j_s\le t}\sum_{C_{j_1}\subset A_{j_1}^N}...\sum_{C_{j_s}\subset A_{j_s}^N}\left(\frac{t}{t-1}\right)^{\sum_{k=1}^sm(C_{j_k})+s}\,\left|\xs\left[\prod_{k=1}^s\prod_{i\in C_{j_k}}\oa\right]\right|\;+...$$
$$
...+\;\sum_{C_1\subset A_1^N}...\sum_{C_t\subset A_t^N}\left(\frac{t}{t-1}\right)^{\sum_{j=1}^tm(C_j)+t}\,\left|\xs\left[\prod_{j=1}^t\prod_{i\in C_j}\oa\right]\right|.\eqno{(3.18)}$$

Let us note that for
$t\ge 3$
$$
\frac{t}{t-1}\;\le\;2 \mbox{\;\;\;and\;\;\;} \left(\frac{t-1}{t-2}\right)^{t-s}\,\left(\frac{t}{t-1}\right)^s\;\le\;4e\;\;(s=0,1,.,t).\eqno{(3.19)}$$

Since the sets
$A_j^N$
$(j=1,2,.,t)$
are mutually disjoint, then the integrands in the last expression are
the distinct products of pairwise distinct functions
$\varphi_i$
$(i=1,2,.,N).$
In addition, the number of these functions is equal to
$m(C_j)$
$(j=1,2,.,t),$
$m(C_{j_1})+m(C_{j_2})$
$(1\le j_1<j_2\le t),..$
$,\sum_{j=1}^tm(C_j)$
for terms of the first,second,..,last sum, respectively.
Therefore the maximal index of them in every integrand is not less than
$m(C_j),$
$m(C_{j_1})+m(C_{j_2}),..$
$\sum_{j=1}^tm(C_j),$
respectively. Finally, by (3.18) and (3.19), there holds
$$
L_N\;\le\;4e\,\left\{\sum_{i=1}^N 2^i\,\sum_{1\le j_1<...<j_s\le i-1}\Bigl|\xs(\varphi_{j_1}...\varphi_{j_s}\varphi_i)\Bigr|\right\},\eqno{(3.20)}$$
where the inner summation is taken over all sets of indices
$1\le j_1<...<j_s\le i-1.$
From this inequality, (3.4), and (3.5) it follows
$$
L_N\;\le\;\frac{4e}{D}\za\yb\;<\;2.$$
So, by (3.15) and (3.16), we get
$$
|I_N|\;\le\;2\oc\;\;\;\;\mbox{for}\;\;\;N=1,2,..$$

Taking into account (3.14), from the last inequality we have
$$
\sh\left(\sum_{i\in A_j}a_i^2\right)^{1/2}\;\le\;6D\oc.$$
Therefore, by Lemma 3.1,
$$
\xy\;\le\;6\sqrt{2}\,D\oc,\eqno{(3.21)}$$
for all positive integers
$t\ge 3$
and
$a=\yz\in l_2.$

Let us prove an analogous inequality in the cases
$t=1$
and
$t=2.$
First, the constructed system
$\qr$
is a Riesz basic sequence. It means that
$$
\Bigl\|\wx\Bigr\|_2\;\asymp\;\|(a_i)\|_2.\eqno{(3.22)}$$
In fact, by assumption,
$\|\varphi_i\|_2=1.$
Therefore (3.5) implies that
$$
\Bigl\|\od a_i\varphi_i\Bigr\|_2^2\;=\;\od a_i^2\,+\,2\sum_{1\le i<j\le m} a_ia_j\xs(\varphi_i\varphi_j)\;\ge$$
$$
\ge\;\od a_i^2\,\left\{\,1\;-\;2\biggl[\sum_{1\le i<j\le m}(\xs(\varphi_i\varphi_j))^2\biggr]^{1/2}\right\}\;\ge\;\frac{1}{2}\|(a_i)\|_2^2,\eqno{(3.23)}$$
for all
$m\in\fg.$
Similarly,
$$
\Bigl\|\od a_i\varphi_i\Bigr\|_2^2\;\le\;\frac{3}{2}\|(a_i)\|_2^2,$$
and (3.22) is proved.

Inequality (3.2) shows that the system
$\{g_i\}_{i=1}^\infty$
is a
$\Lambda(\infty)$
-system.
Therefore
$\qr$
hase the same property, and hence it is a Banach system [3,Cor.1.3.1].
It means that
$$
\|a\|_2\;\le\;M\,\Bigl\|\od a_i\varphi_i\Bigr\|_1\;\;\;\;\mbox{for}\;\;\;a=\yz\in l_2,\eqno{(3.24)}$$
where a constant
$M>0$
depends only on
$D.$

The inequality
$\|a\|_2\,\le\,\|a\|_1$
implies that
$\gh(1,a;l_1,l_2)\,=\,\|a\|_2.$
Hence, using the properties of the
$\gh$
-functional, (3.23), and (3.24) we obtain:
$$
\gh_{1,2}(\sqrt{2},a)\;\le\;\sqrt{2}\,\gh_{1,2}(1,a)\;\le\;2\|\varphi\|_2,\eqno{(3.25)}$$
and
$$
\gh_{1,2}(1,a)\;\le\;M\|\varphi\|_1.\eqno{(3.26)}$$

Suppose now that
$t\ge 1$
is an arbitrary real number. Choose a positive integer
$t_0$
such that
$t_0\le t<t_0+1.$
Then (3.21), (3.25), and (3.26) yield
$$
\xy\;\le\;\sqrt{t/t_0}\,\gh_{1,2}(\sqrt{t_0},a)\;\le\;\sqrt{2}\,\gh_{1,2}(\sqrt{t_0},a)\;\le\;C\oc,$$
where
$C=\,\sqrt{2}\max(2,M,6\sqrt{2}\,D)$
depends only on
$D.$

To conclude the proof, it remains to note that (3.3) holds,
in particular, for the system
$\qr.$
\end{proof}

\begin{theor}
A system
$\bc$
of r.v. defined on a probability space
$\ab$
contains a subsystem
$\qr$
equivalent in distribution to the Rademacher system
on the segment
$[0,1]$
if and only if there exists a subsystem
$\cd\subset\mn$
that satisfies the following conditions:

1) $|f_{n_k}(\omega)|\le D$
$(k=1,2,..;$
$\de$);

2) $f_{n_k}\to 0$
weakly in
$L_2;$

3) $d=\,\inf_{k=1,2,..}\|f_{n_k}\|_2>0.$

Moreover, a constant of the equivalence
$\xc\stackrel{\vr}{\sim}\pq$
depends only on
$D$
and
$d.$
\end{theor}
\begin{proof}
First suppose there exists a subsystem
$\cd\subset\mn$
satisfying conditions 1) --- 3).
By Theorem 3.2, we can select a subsystem
$\qr\subset\bc$
such that equivalence (3.1) holds. Using Theorem 2.1, we obtain:
$\xc\stackrel{\vr}{\sim}\pq.$

Conversely, suppose that
$\qr\subset\bc$
and
$\xc\stackrel{\vr}{\sim}\pq.$
Clearly, the Rademacher system possess the
properties similar to 1) --- 3). Therefore
$\qr,$
also, satisfies 1) and 3). As we remarked in Section 1,
the Rademacher system is a
$\Lambda(\infty)$
-system. Hence
$\qr$
have this property, too. The possibility of the selection of subsystem
$\cd\subset\qr$
satisfying condition 2) is a consequence of well-known Stechkin's theorem
[3,Th.1.3.1]. This completes the proof.
\end{proof}

The following result is an immediate consequence of Theorem 3.3.
\begin{theor}
Suppose
$\bc$
is an orthonormal sequence of r.v. defined on a probability space
$\ab,$
$|f_n(\omega)|\le D$
($n=1,2,..;\de).$
Then it contains a subsequence
$\qr$
equivalent in distribution to the Rademacher system. A constant
of this equivalence depends only on
$D.$
\end{theor}
\begin{rem}
In [18,${\cal x}\,3.2$], G.Alexits posed the problem if any total orthogonal
system of functions contains an infinite multiplicative or even strongly
multiplicative subsystem. Theorem 3.4 give a positive solution
of a problem similar to it: any uniformly bounded orthogonal system of
r.v. contains a subsystem equivalent in distribution to the Rademacher
system (consisting of independent functions).
\end{rem}

As already stated in Section 1, if
$\xc\stackrel{\vr}{\sim}\pq,$
then
$\xc$
is a Sidon system. It turns out a such system possess, in a definit sense,
a strengthened Sidon property.
\begin{cor}
Suppose that a sequence
$\bc$
of r.v. contains a subsequence
 $\cd,$
satisfying conditions 1) --- 3) of Theorem 3.3. Then
there are a subsystem
$\qr\subset\bc$
and constants
$\alpha_1,\alpha_2,\alpha_3>0$
depending only on
$D$
and
$d$
with following property: for an arbitrary polynomial
$$
T(\omega)=\;\od a_i\varphi_i(\omega)$$
there exists a set
$E=E(T)\subset\Omega$
such that
$\vr(E)>\alpha_12^{-m}$
and
$$
|T(\omega)|\;\ge\;\alpha_2\,\|T\|_\infty\;\ge\;\alpha_3\,\od |a_i|\;\;\;\;\mbox{for}\;\;\;\omega\in E.$$
\end{cor}
\begin{proof}
By Theorem 3.3,
$\mn$
contains a subsequence
$\xc$
that is equivalent in distribution to the Rademacher system. It means that
$$
C^{-1}|\{|\tilde{T}(x)|>Cz\}|\;\le\;\vr\{|T(\omega)|>z\}\;\le\;C|\{|\tilde{T}(x)|>C^{-1}z\}|,$$
where
$\tilde{T}(x)=\;\od a_ir_i(x)$
$(x\in [0,1])$
and
$C>0$
is a constant independent of
$z>0.$

In particular, it yields
$$
C^{-1}\|T\|_\infty\;\le\;\|\tilde{T}\|_\infty\;\le\;C\|T\|_\infty.$$
Therefore, by definition of Rademacher functions,
\begin{multline*}
\vr\{|T(\omega)|>(2C^2)^{-1}\|T\|_\infty\}\;\ge\;\vr\{|T(\omega)|>(2C)^{-1}\|\tilde{T}\|_\infty\}\;\ge\\
\ge\;C^{-1}\,|\{|\tilde{T}(x)|>2^{-1}\|\tilde{T}\|_\infty\}|\;\ge\;C^{-1}2^{-m}.
\end{multline*}

Denote
$\alpha_1=\,C^{-1},$
$\alpha_2=\,(2C^2)^{-1},$
and
$E\,=\,E(T)=\;\{\de:\:|T(\omega)|>\alpha_2\|T\|_\infty\}.$
Then if
$\omega\in E,$
$$
|T(\omega)|\;>\;\alpha_2\|T\|_\infty\;\ge\;\alpha_2C^{-1}\|\tilde{T}\|_\infty\;=\;\alpha_3\od |a_i|,$$
where
$\alpha_3=\,\alpha_2C^{-1}.$
\end{proof}
\vskip 0.2cm

Remind that there is a local version of the Sidon property. Namely,
a sequence
$\bc$
of functions measurable on the segment
$[0,1]$
is called a Sidon-Zygmund system if there exists a set
$E\subset [0,1]$
such that
$|E|>0$
and
the condition
$\op a_nf_n(x)\in L_\infty(I)$
implies
$\ik\in l_1,$
for any segment
$I$
such that
$|I\cap E|>0.$
We say that
$\bc$
is a Sidon-Zygmund system in the restricted sense if, in addition,
we can take
$E=\,[0,1].$

The last results of this section strengthen the theorem of Gaposhkin
about the selection of subsystems with the Sidon-Zygmund property (in the usual and
the restrictive sense) [3,Th.1.4.2]. By
$\chi_E$
we denote here the indicator function of a set
$E\subset [0,1],$
i.e.,
$\chi_E(x)=\,1$
for
$x\in E$
and
$\chi_I(x)=\,0$
for
$x\not\in E.$
\begin{theor}
Suppose a sequence
$\bc$
of functions measurable on the segment
$[0,1]$
contains a subsequence
$\cd$
satisfying conditions 1) --- 3) of Theorem 3.3. There exist a subsystem
$\xc\subset\mn$
and a set
$E\subset [0,1],$
$|E|>0$
having the following property:

for any segment
$I\subset [0,1]$
such that
$|I\cap E|>0$
and for some positive integer
$k_0=k_0(I)$
depending on
$I$
the sequence
$\{\varphi_i\chi_I\}_{i=k_0}^\infty$
is equivalent in distribution to the Rademacher system on
$[0,1].$
A constant of this equivalence depends on
$D,$
$d,$
and
$I.$
\end{theor}
\begin{proof}
By Theorem 2.1, it is enough to find a subsystem
$\xc\subset\mn$
and a set
$E\subset [0,1],$
$|E|>0$
such that for every segment
$I\subset [0,1],$
$|I\cap E|>0$
and for some
$k_0=k_0(I)$
we have
$$
\Bigl\|\za a_i\varphi_{k_0+i-1}\chi_I\Bigr\|_t\;\asymp\;\xy\;\;\;\mbox{for}\;\;\;a=\yz\;\;\;\mbox{and}\;\;\;t\ge 1.$$
(a constant of this equivalence depends on
$D,$
$d,$
and
$I$).

Next, we can use arguments similar to those in the proof of Theorem 3.2.
Hence we make only some remarks.

First of all, it is clear that inequality (3.2) holds for functions
$g_i\chi_I$
$(I$
is any segment from
$[0,1]).$
Next, as in the proof of Theorem 1.4.2 from [3], we can select
a subsequence
$\xc\subset\{g_i\}$
such that the relations similar to (3.5) hold for the integrals
taken over an arbitrary segment
$I\subset [0,1].$
Moreover, by well-known Marcinkiewicz's lemma (see [19[ or [3,L.1.2.5]),
we can assume that there is a some set
$E\subset [0,1]$
such that
$|E|>0$
and for any subset
$F\subset E,$
$|F|>0$
there holds:
$$
\liminf_{i\to\infty}\int_F\varphi_i^2(x)\,dx\;>\;0.\eqno{(3.27)}$$
So, we can repeat next the arguments of the proof of Theorem 3.2  replacing
only the integrals over
$[0,1]$
with the integrals over a segment
$I$
such that
$|I\cap E|>0.$
\end{proof}
\begin{theor}
Suppose a sequence
$\bc$
of functions measurable on
$[0,1]$
contains a subsequence
$\cd$
satisfying conditions 1) and 2) of Theorem 3.3. Besides, if
$F\subset[0,1]$
such that
$|F|>0,$
then
$$
\liminf_{k\to\infty}\int_F f_{n_k}^2(x)\,dx\;>\;0\eqno{(3.28)}$$

Then there exists a subsystem
$\xc\subset\mn$
with the following property: for any segment
$I\subset [0,1],$
$|I|>0$
can be found a positive integer
$k_0=k_0(I)$
such that the sequence
$\{\varphi_i\chi_I\}_{i=k_0}^\infty$
is equivalent in distribution to the Rademacher system on
$[0,1].$
\end{theor}
\begin{proof}
in view of (3.28) we can assume that the consructed subsystem
$\qr$
satisfies (3.27) for an arbitrary set
$F\subset [0,1]$
of positive Lebesgue measure. Next, we argue as in the proofs of Theorems 3.7
and 3.2.
\end{proof}
\vskip 0.6cm

{\bf 4. Density of subsystems equivalent in distribution to
the Rademacher system.}
\setcounter{section}{4}
\setcounter{theor}{0}
Now we consider a problem of selection of lacunary subsets from finite
orthogonal sets
$\abf$
of functions defined on the segment
$[0,1].$
We prove that a such set contains a subset
$\abe$
of "logarithmic density" (i.e.,
$s\ge C\log_2N)$
that is equivalent in distribution to the set of first
$s$
Rademacher functions with a constant independent of
$N.$
This strengthens Kashin's theorem on selection of Sidon subsets [4].
Let us note that the analogous
$\Lambda(p)$
-problem was solved by J.Bourgain [20].
\begin{theor}
Let
$\abf$
be an orthonormal set of functions defined on the segment
$[0,1],$
$|f_n(x)|\le D$
$(n=1,2,.,N).$
There exist a subset
$\abe$
$(1\le n_1<n_2<...<n_s\le N)$
and a constant
$C>0$
depending only on
$D$
such that
$s\ge\max\{[1/6\log_2N],1\}$
and
$$
C^{-1}\Bigl|\Bigl\{\Bigl|\sum_{i=1}^{s} a_i r_i(x)\Bigr|>Cz\Bigr\}\Bigr|\;\le\;\Bigl|\Bigl\{\Bigl|\sum_{i=1}^{s} a_i \abd(x)\Bigr|>z\Bigr\}\Bigr|\;\le\;C\Bigl|\Bigl\{\Bigl|\sum_{i=1}^{s} a_i r_i(x)\Bigr|>C^{-1}z\Bigr\}\Bigr|,$$
for all
$a_i\in\hj$
$(i=1,2,.,s)$
and
$z>0.$
\end{theor}

We need two lemmas; the first was proved by B.S.Kashin [17,L.8.4.1].

\begin{lemma}
Suppose a set of functions
$\abf$
$(\log_2N\ge 6)$
satisfies the conditions of Theorem 4.1. Then there exists a subset
$\abe\subset\abf$
$(1\le n_1<n_2<...<n_s\le N)$
such that
$s\ge [1/6\log_2N]$
and
$$
\cdh\left\{\xs\left[\abk\biggl(\frac{\abd(x)}{D}\biggr)^\abi\right]\right\}^2\;\le\;10^{-s},\eqno{(4.1)}$$
where
$\abg$
is the collection of all
$\theta=(\abi)_{i=1}^s$
such that
$$
a)\;\abi=\,0,1 \hbox{\;\; or\;\;} 2\;\;;\;\;b)\;\sum_{i:\,\abi=1} 1\,\ge\,1\;\;;\;\;c)\;\sum_{i:\,\abi=2} 1\,\le\,1.$$
\end{lemma}

The second lemma allows, under certain assumptions, to extend an
uniformly bounded set of functions from the segment
$[0,1]$
to some more large segment, so that the constructed set is
uniformly bounded and multiplicative.
Analogous results about extending to an orthogonal set see in
[17,${\cal x}\,7.1$].

\begin{lemma}
Let
$\{g_i\}_{i=1}^s$
be a set of functions on
$[0,1],$
$|g_i(x)|\le D$
$(i=1,2,.,s;\,x\in [0,1]),$
and
$$
\max_{\theta\in\abg'}\biggl|\xs\biggl[\abk\biggl(\frac{g_i(x)}{D}\biggr)^\abi\biggr]\biggr|\;<\;2^{-s},\eqno{(4.2)}$$
where
$\abg'$
is the collection of all
$\theta=\,(\abi)_{i=1}^s\in\abg$
such that
$\abi=0$
or
$1.$

Then there exists a multiplicative set
$\abj$
of functions defined on the segment
$[0,2]$
such that
$h_i(x)\chi_{[0,1]}(x)=g_i(x)$
and
$\|h_i\|_{L_\infty[0,2]}\le D$
$(i=1,2,.,s).$
\end{lemma}
\begin{proof}
For any positive integer
$s,$
we have
$$
[1,2)=\,\bigcup_{k=1}^{2^s}\abl\,,\;\;\;\mbox{where}\;\;\;\abl=[a_{k-1},a_k)\,,\,a_k=1+k2^{-s}\,(k=1,2,.,2^s).$$
Let us fix an one-to-one mapping from the collection
$\abg'$
to the set of intervals
$\{\abl\}_{k=1}^{2^s-1}.$

Suppose
$\theta=(\theta_i)_{i=1}^s$
is the element of the collection
$\abg'$
corresponding to an interval
$\abl.$
Denote by
$\{i_j\}_{j=1}^m$
$(1\le i_1<i_2<...<i_m\le s)$
the set of all
$i=1,2,.,s$
such that
$\abi=1.$

For every
$\alpha\in\abl,$
we define
$$
u_\alpha(x)=\;\chi_{[a_{k-1},\alpha)}(x)\,-\,\chi_{[\alpha,a_k]}(x)\;\;\;\;\mbox{for}\;\;\;x\in\abl.\eqno{(4.3)}$$
Then the image of the function
$$
v(\alpha)=\;\int_{\abl}u_\alpha(x)\,dx\;=\;2\alpha\,-\,a_{k-1}\,-\,a_k\;\;\;\;\mbox{for}\;\;\;\alpha\in\abl$$
includes the segment
$[-2^{-s},2^{-s}].$
By condition (4.2), the functions
$\abm_i(x)=\,g_i(x)/D$
satisfy the inequality:
$|\xs(\abm_{i_1}\abm_{i_2}...\abm_{i_m})|\,<\,2^{-s}.$
Therefore it can be found
$\abn\in (a_{k-1},a_k)$
such that
$$
v(\abn)=\;-\xs(\abm_{i_1}\abm_{i_2}...\abm_{i_m}).\eqno{(4.4)}$$

First suppose
$m\ge 2.$
Denote by
$d_j(x)$
$(j=1,2,.,m-1)$
step functions defined for
$x\in\abp$
and having the following properties:
$$
1)\;|d_j(x)|\,\equiv\,1\;\;;\;\;2)\;\int_{a_{k-1}}^\abn d_{j_1}d_{j_2}...d_{j_l}\,dx\,=\,0\;(1\le j_1<j_2<...<j_l\le m-1)\eqno{(4.5)}$$
(for example, a such functions can be constructed using the Rademacher
functions). Let
$c_j(x)$
$(j=1,2,.,m-1)$
be an analogous system, but defined on the interval
$\abq.$

Let us define the functions
$h_i(x)$
as follows.
If
$x\in\abl$
$(k=1,..,2^s-1),$
then
$$
h_i(x)\,\equiv\,0\,(i\ne i_j,j=1,2,.,m)\,,\,h_{i_j}(x)\,=$$
$$
=\,D\,[d_j(x)\chi_\abp(x)\,+\,c_j(x)\chi_\abq(x)]\,(j=1,2,.,m-1),\eqno{(4.6)}$$
and
$$
h_{i_m}(x)\,=\,D\,\Bigl[\prod_{j=1}^{m-1}d_j(x)\chi_\abp(x)\,-\,\prod_{j=1}^{m-1}c_j(x)\chi_\abq(x)\Bigr].\eqno{(4.7)}$$
By (4.5) and (4.4), we have
$$
\int_\abl h_{i_1}h_{i_2}...h_{i_m}\,dx\;=D^m\int_\abl u_{\abn}(x)\,dx\;=\;D^mv(\abn)\;=\;-\xs(g_{i_1}g_{i_2}...g_{i_m}).\eqno{(4.8)}$$
At the same time, from the definition of functions
$h_i$
$(i=1,2,.,s)$
and (4.5) it follows
$$
\int_\abl h_{i_1'}h_{i_2'}...h_{i_l'}\,dx\;=\;0,\eqno{(4.9)}$$
for every another set of indices
$\{i_j'\}_{j=1}^l$
$(1\le i_1'<...<i_l'\le s).$

In the case
$m=1$
we define
$$
h_i(x)\,\equiv\,0\;\;\;\;\mbox{if}\;\;\;i\ne i_1\;\;\;\;\mbox{and}\;\;\;h_{i_1}(x)\,=\,Du_\abn(x).\eqno{(4.10)}$$
Then equalities (4.8) and (4.9) hold again.

So, the functions
$h_i(x)$
are defined on
$\cup_{k=1}^{2^s-1}\Delta_k.$
In addition, let
$h_i(x)=\,g_i(x)$
for
$x\in [0,1]$
and
$h_i(x)=\,0$
for
$x\in\Delta_{2^s}.$
Since (4.3) and (4.5) --- (4.10) yield that
$\abj$
is a multiplicative set on the segment
$[0,2]$
and
$|h_i(x)|\le D,$
then this completes the proof.
\end{proof}

\begin{proof}[Proof of Theorem 4.1]
In view of Theorem 2.1 and Remark 2.2 it is sufficient to prove
that there exists a subset
$\abe$
such that
$s\ge\max\{[1/6\log_2N],1\}$
and
$$
\Bigl\|\abc a_i\abd\Bigr\|_t\;\asymp\;\xy\;\;(t\ge 1\,,\,a=\abr),\eqno{(4.11)}$$
with a constant depending only on
$D.$
As before, here
$\xy=\,\gh(\sqrt{t},a;l_1,l_2).$

Using Lemma 4.2, we can select a subset
$\abe$
$(s\ge\max\{[1/6\log_2N],1\})$
satisfying condition (4.1). Since
$\abg'\subset\abg,$
then (4.2) holds also for functions
$g_i=\,f_{n_i}.$
Hense, by Lemma 4.3, the set
$\abe$
can be extended to a multiplicative set
$\abj$
on the segment
$[0,2],$
$|h_i(x)|\le D.$
It is easy to check that the set of functions
$h_i'(x)=\,h_i(2x)$
$(i=1,..,s)$
is a multiplicative system on the segment
$[0,1]$
and
$|h_i'(x)|\le D.$
Therefore, by Corollary 3 from [10],
$$
\Bigl\|\abc a_ih_i'\Bigr\|_{L_t[0,1]}\;\le\;C_1\sqrt{t}\|a\|_2,$$
where
$t\ge 1,$
$a=\abr,$
and a positive constant
$C_1$
depends only on
$D.$
Hence,
$$
\Bigl\|\abc a_i\abd\Bigr\|_t\;\le\;2C_1\sqrt{t}\|a\|_2\eqno{(4.12)}$$
Using the last inequality and arguing as in the proof of Theorem 3.2
(see inequality (3.20)), we obtain:
$$
\Bigl\|\abc a_i\abd\Bigr\|_t\;\le\;C\xy\;\;(t\ge 1),\eqno{(4.13)}$$
where
$C=\,\max(D,2C_1)$
depends only on
$D.$

Let us prove the opposite inequality. Denote
$$
a=\,\abr  \hbox{\;\;\; and\;\;\;   }  P(x)=\,\abc a_i\abd(x).$$

At first, suppose that
$t\in\fg,$
$t\ge 3.$
For
$s\le 16D^2,$
by the definition of the
$\gh$
-functional, we have:
$$
\|P\|_t\;\ge\;\|P\|_2\;=\;\|a\|_2\;\ge\;\frac{1}{\sqrt{s}}\abc |a_i|\;\ge\;\frac{1}{4D}\xy.\eqno{(4.14)}$$
Therefore it is sufficient to consider the case when
$$
s\;>\; 16D^2.\eqno{(4.15)}$$
Let
$\kl$
be an arbitrary partition of the set
$\{1,2,.,s\}$
(probably, some of them are empty) and
$$
\abt b_i^2\;\le\;1\;\;\;\;\mbox{for}\;\;\;j=1,2,.,t.\eqno{(4.16)}$$
Introduce the Riesz products
$$
R_s(x)=\;\abk\biggl(1\,+\,\frac{b_i}{D}\abd(x)\biggr)$$
and the integral
\begin{multline*}
I_s\;=\;\int_0^1P(x)R_s(x)\,dx\;=\;\abc a_i\,\int_0^1\abd(x)\,dx\;+\\
+\;\abc a_i\,\int_0^1\abd(x)\abh\abv\biggl[\frac{b_k}{D}\abu(x)\biggr]^\abw\,dx.
\end{multline*}
After regrouping its terms we get
$$
I_s\;=\;D^{-1}\abc a_ib_i\;+\;\abc a_i\abx,\eqno{(4.17)}$$
where
$$
\abx=\;\int_0^1\abd(x)\,dx\,+\,\int_0^1\abd(x)\sum_{\theta\in\abg'(i)}\abv\biggl[\frac{b_k}{D}\abu(x)\biggr]^\abw\,dx,$$
$$
\abg'(i)=\;\abg'\setminus \{\theta^i\}\;\;\;\mbox{where}\;\;\;\theta^i=(\theta_j^i),\theta_i^i=1,\theta_j^i=0(j\ne i).$$
Since
$|b_i|\le 1$
$(i=1,2,.,s),$
then (4.1) and the Cauchy-Bunyakovskii inequality yield
$$
\abc|\abx|\;\le\;D\cdh\biggl|\int_0^1\abv\biggl[\frac{\abu(x)}{D}\biggr]^\abw\,dx\biggr|\;\le\;D(s2^s)^{1/2}10^{-s/2}\;=\;D(s5^{-s})^{1/2}.$$
By (4.15),
$(s5^{-s})^{1/2}\,\le\,8/s\,<\,2^{-1}D^{-2},$
and we have
$$
\abc|\abx|\;\le\;\frac{1}{2D}.\eqno{(4.18)}$$

Let
$b_i$
$(i\in A_j,\,j=1,..,t)$
satisfy (4.16) and
$$
\abt a_ib_i\;=\;\Bigl(\abt a_i^2\Bigr)^{1/2}.$$
Then from (4.17) and (4.18) we obtain
\begin{multline*}
I_s\;\ge\;\frac{1}{D}\sh\Bigl(\abt a_i^2\Bigr)^{1/2}\,-\,\sh\Bigl|\abt a_i\abx\Bigr|\;\ge\;\frac{1}{D}\sh\Bigl(\abt a_i^2\Bigr)^{1/2}\,-\\
-\,\sh\Bigl(\abt a_i^2\Bigr)^{1/2}\,\Bigl(\abt |\abx|\Bigr)\;\ge\;\frac{1}{2D}\sh\Bigl(\abt a_i^2\Bigr)^{1/2}.
\end{multline*}
Thus, by (2.2), we have
$$
I_s\;\ge\;\frac{1}{2\sqrt{2}D}\,\xy\;\;\;\;\mbox{for positive integer}\;\;\;t\ge 3.\eqno{(4.19)}$$

Now, let us estimate the integral
$I_s$
from above. By H\"older's inequality,
$$
|I_s|\;\le\;\|R_s\|_{t'}\|P\|_t,$$
where
$t'=\,t/(t-1).$

Arguing as in the proof of Theorem 3.2, we can prove that
$$
\|R_s\|_{t'}\;\le\;4e2^s\,\abh\biggl|\xs\biggl\{\abk\biggl[\frac{\abd(x)}{D}\biggr]^\abi\bigg\}\biggr|\;\;\;\;\mbox{for}\;\;\;t\ge 3.$$
So, inequality (4.1), the embedding
$\abg'\subset\abg,$
and the Cauchy-Bunyakovskii inequality imply
$\|R_s\|_{t'}\,\le\,4e.$
Therefore
$|I_s|\,\le\,4e\|P\|_t,$
and from (4.19) it follows
$$
\xy\;\le\;8\sqrt{2}eD\|P\|_t\;\;\;\;\mbox{for positive integers}\;\;\;t\ge 3.\eqno{(4.20)}$$
Analogous relations hold also for
$t=2$
and
$t=1.$
In fact, since
$\abf$
is an orthonormal set and
$\gh_{1,2}(1,a)=\|a\|_2,$
then
$$
\gh_{1,2}(\sqrt{2},a)\;\le\;\sqrt{2}\gh_{1,2}(1,a)=\sqrt{2}\|P\|_2.\eqno{(4.21)}$$
In addition, as in the proof of Theorem 3.2, using (4.12) and Corollary 1.3.1
from [3], we obtain
$$
\gh_{1,2}(1,a)\;=\;\|a\|_2\;\le\;M\|P\|_1,\eqno{(4.22)}$$
where a constant
$M>0$
depends only on
$D.$

Finally, we extend, as usual, (4.20) --- (4.22) to all real
$t\ge 1$
and get
$$
\xy\;\le\;C\|P\|_t,$$
where
$C=\,\sqrt{2}\max(M,\sqrt{2},8\sqrt{2}eD).$
Note that this constant depends only on
$D.$
Therefore in view of (4.13) equivalence (4.11) holds with a conctant
depending only on
$D.$
This completes the proof.
\end{proof}
\begin{rem}
The "logarithmic" order conditions of the density of subsystem
$\abe$
are exact. Indeed, in the paper [21] S.B.Stechkin proved that
$$
\sum_{k:n_k<N}\,1\;\le\;C\ln N\;\;\;\;\mbox{for}\;\;\;N=2,3,..$$
whenever a sequence
$\{\sqrt{2}\cos 2\pi n_kx\}_{k=1}^\infty\;(x\in [0,1])$
is a Sidon system.
\end{rem}
\begin{rem}
Lacunary subsystems equivalent in distribution to the Rademacher system
(to be more precise, distributions of corresponding polynomials)
are best possible. Indeed, it is readily seen that
$$
\Bigl|\Bigl\{\Bigl|\od a_ir_{n_i}(x)\Bigr|>z\Bigr\}\Bigr|\;=\;\Bigl|\Bigl\{\Bigl|\od a_ir_i(x)\Bigr|>z\Bigr\}\Bigr|,$$
for any subsystem
$\{r_{n_i}\}_{i=1}^\infty\subset\{r_n\}_{n=1}^\infty,$
$m\in\fg,$
$a_i\in\hj$
$(i=1,2,.,m),$
and
$z>0.$
\end{rem}
\vskip 0.6cm

{\bf 5. $\gh$
-closed representability of Banach couples.}
\setcounter{section}{5}
\setcounter{theor}{0}
At first, we introduce the following concept. We shall say that a Banach
couple
$(X_0,X_1)$
is
$\gh$
-closed representable in a Banach couple
$(Y_0,Y_1)$
if there exists a linear operator
$T:\,X_0+X_1\to Y_0+Y_1$
such that

$a)$
$T$
is a injektive bounded operator from
$X_0$
into
$Y_0$
and from
$X_1$
into
$Y_1;$

$b)$ with a constant
$C>0$
independent of all
$x\in X_0+X_1$
and
$t>0$
$$
\gh(t,Tx;Y_0,Y_1)\;\asymp\;\gh(t,x;X_0,X_1).$$

As was proved in the paper [22] (see also [23]), the couple
$(l_1,l_2)$
is
$\gh$
-closed representable in the couple of function spaces
$(L_{\infty},G).$
By
$G$
we denote here the closure of space
$L_{\infty}$
in the Orlicz space
$L_N,$
$N(u)=\,e^{u^2}-1.$
Moreover, we can take as an appropriate operator 
$$
Ta(x)=\;\sum_{n=1}^{\infty}a_nr_n(x)\;\;\;\;\mbox{for}\;\;\;a=(a_n)_{n=1}^{\infty}\,\in\,l_2.$$

Suppose that a sequence of functions
$\bc$
measurable on
$[0,1]$
satisfies conditions 1) --- 3) of Theorem 3.3. Then, this sequence
contains a subsequence
$\xc$
that is equivalent in distribution to the Rademacher system. In the paper
[22], estimates of the tail distribution of Rademacher polynomials are
used, only. For this reason all results of [22]
hold also for
$\xc.$
In particular, every such operator
$$
T_{\varphi}a(x)=\;\wx(x)\;\;(a=\yz\in l_2)$$
realizes the
$\gh$
-closed representability of the Banach couple
$(l_1,l_2)$
in the couple
$(L_\infty,G).$

Using the results of Section 3, we shall prove here a negative
assertion closely connected to the previous example.
\begin{theor}
The Banach couple
$(l_1,l_2)$
is not
$\gh$
-closed representable in a Banach couple
$(L_\infty,L_2)$
of spaces of r.v. defined on a probability space
$\ab.$
\end{theor}
\begin{proof}
Assune the converse. Let
$T:\,l_2\to L_2$
be a linear operator such that
$$
\|Ta\|_{\infty}\,\asymp\,\uv\;\;\;\mbox{for}\;\;a\in l_1,\;\;\; \|Ta\|_2\,\asymp\,\|a\|_2\;\;\;\mbox{for}\;\;a\in l_2,\eqno{(5.1)}$$
and
$$
C_1^{-1}\gh_{1,2}(t,a)\;\le\;\gh(t,Ta;L_\infty,L_2)\;\le\;C_1\gh_{1,2}(t,a)\;\;\;\mbox{for}\;\;a\in l_2\;\;\;\mbox{and}\;\;t>0.\eqno{(5.2)}$$
Here, as before,
$\gh_{1,2}(t,a)=\,\gh(t,a;l_1,l_2).$

Denote
$f_n=\,Te_n,$
for
$e_n=\,(\delta_n^j),$
where
$\delta_n^n=1$
and
$\delta_n^j=0,$
if
$j\ne n.$
Using (5.1), we obtain:
$$
D^{-1}\uv\;\le\;\Bigl\|\op a_nf_n\Bigr\|_\infty\;\le\;D\uv\;\;\;\;\mbox{if}\;\;\;a\in l_1$$
and
$$
d\|a\|_2\;\le\;\Bigl\|\op a_nf_n\Bigr\|_2\;\le\;d^{-1}\|a\|_2\;\;\;\;\mbox{if}\;\;\;a\in l_2,\eqno{(5.3)}$$
for some constants
$D>0$
and
$d>0.$
Hence, in particular,
$$
|f_n(\omega)|\;\le\;D\;\;(n=1,2,..;\de)\;\; \mbox{\;\;\;   and\;\;\;       }\;\; \|f_n\|_2\;\ge\;d.\eqno{(5.4)}$$

In addition, from (5.3) it follows that
$f_n\ne 0$
$(n=1,2,..)$
and
$$
\Bigl\|\sum_{n=1}^k a_nf_n\Bigr\|_2\;\le\;d^{-2}\Bigl\|\ef\Bigr\|_2,$$
for any
$k,m\in\fg,$
$k\le m.$
So,
$\bc$
is a basis in the closed linear span
$Z$
generated by this system in the space
$L_2$
[26,Pr.1.a.3]. In other words,
$\bc$
is a Riesz basic sequence, and therefore 
$\xs(gf_n)\to 0$
as
$n\to\infty,$
for every
$g\in Z$
[3,${\cal x}\,1.1$].

For an arbitrary
$g\in L_2.$
we write:
$g=\,g_1+g_2,$
$g_1\in Z,$
$g_2\in Z^{\perp},$
where
$Z^{\perp}$
is the orthogonal complement to
$Z.$
Then, obvious,
$\xs(gf_n)=\xs(g_1f_n)\to 0$
as
$n\to\infty.$
Finally, we have
$$
f_n\to 0\;\; \mbox{\;    weakly\;\; in  }\;\; L_2\;\;\;\;\mbox{as}\;\;n\to\infty.\eqno{(5.5)}$$

Relations (5.4) and (5.5) show that the system
$\bc$
satisfies conditions 1) --- 3) of Theorem 3.3. Therefore there
exists a subsystem
$\xc\subset\mn$
that is equivalent in distribution to the Rademacher system on the segment
$[0,1].$
This and Khintchine's inequality (see [25] or [26,${\cal x}\,4.5$])
imply
$$
C_2^{-1}\|a\|_2\;\le\;\Bigl\|\wx\Bigr\|_p\;\le\;C_2\|a\|_2\;\;\;\;\mbox{for}\;\;\;a=\yz\in l_2,\eqno{(5.6)}$$
where a positive constant
$C_2$
depends only on
$D,$
$d,$
and
$p\in [1,\infty).$

At the same time, (5.2) yields
$$
C_1^{-1}\gh_{1,2}(t,a)\;\le\;\gh(t,\wx;L_\infty,L_2)\;\le\;C_1\gh_{1,2}(t,a)\;\;\;\mbox{for}\;\;a\in l_2\;\;\mbox{and}\;\;t>0.$$
Let
$0<\theta<1$
and
$p=\,2/\theta.$
By the last inequalities and
[6,Th.5.2.1], the application of the real interpolation method
$(\cdot,\cdot)_{\theta,p}$
(see Introduction) to the Banach couples
$(L_\infty,L_2)$
and
$(l_1,l_2),$
gives
$$
\Bigl\|\wx\Bigr\|_p\;\asymp\;\|a\|_{r,p}=\;\left\{\za(a_i^*)^pi^{p/r-1}\right\}^{1/p}.$$
Here
$r=\,2/(2-\theta)\,<2,$
and a constant of this equivalence depends only on
$C_1$
and
$\theta.$
Since
$\|a\|_{r,p}\not\asymp\|a\|_2,$
then the last contradicts with inequalities (5.6), if
$p=\,2/\theta.$
This completes the proof.
\end{proof}

The author is grateful to professor B.S.Kashin for the posing of the
problem considered in Section 4.

\end{document}